\documentclass{article}
\usepackage{spconf,amsmath,graphicx}

\usepackage{bbm}

\usepackage{graphicx}%
\usepackage{multirow}%
\usepackage{amsmath,amssymb,amsfonts}%
\usepackage{amsthm}%
\usepackage{mathrsfs}%
\usepackage{xcolor}%
\usepackage{textcomp}%
\usepackage{manyfoot}%
\usepackage{booktabs}%
\usepackage{algorithm}%
\usepackage{algorithmicx}%
\usepackage{algpseudocode}%
\usepackage[left=2cm, right=2cm, top=2.5cm, bottom=3cm]{geometry}

\usepackage{threeparttable}

\usepackage{cite}
\usepackage{epstopdf}%

\usepackage{here}
\usepackage{bm}
\usepackage{algorithm}
\usepackage{algpseudocode}
\usepackage{comment}
\usepackage{cases}
\usepackage{mathtools}
\usepackage{enumitem}
\usepackage{mathabx}
\usepackage{stmaryrd}
\usepackage{tikz}
\usetikzlibrary{arrows.meta}
\usetikzlibrary{positioning}
\usetikzlibrary{decorations.pathreplacing,calc}
\usepackage{csvsimple}
\usepackage{lscape}

\usepackage{etoolbox}
\apptocmd{\liminf}{\limits}{}{}
\apptocmd{\limsup}{\limits}{}{}

\usepackage{tabularx}
\usepackage{booktabs}

\newcommand{\argmin}{\mathop{\mathrm{argmin}}\limits}

\newcommand{\inprod}[2]{{\left\langle #1,#2 \right\rangle}}

\newcommand{\TT}{\mathsf{T}}

\newcommand{\prox}[1]{\mathrm{prox}_{#1}}

\newcommand{\moreau}[2]{{}^{#2}#1}

\newcommand{\exR}{\mathbb{R}\cup \{+\infty\}}

\newcommand{\dom}[1]{\mathrm{dom}(#1)}

\newcommand{\norm}[1]{\left\lVert #1 \right\rVert}

\newcommand{\Fsubdiff}{\partial_{F}}

\newcommand{\dist}{\mathop{\mathrm{dist}}}

\makeatletter
\newcommand{\doublewidetilde}[1]{{%
  \mathpalette\double@widetilde{#1}%
}}
\newcommand{\double@widetilde}[2]{%
  \sbox\z@{$\m@th#1\widetilde{#2}$}%
  \ht\z@=.9\ht\z@
  \widetilde{\box\z@}%
}
\makeatother

\theoremstyle{plain}%
\newtheorem{theorem}{Theorem}[section]
\newtheorem{lemma}[theorem]{Lemma}

\theoremstyle{definition}
\newtheorem{definition}[theorem]{Definition}
\newtheorem{fact}[theorem]{Fact}

\newtheorem{problem}[theorem]{Problem}
\newtheorem{remark}[theorem]{Remark}

\theoremstyle{remark}

\DeclarePairedDelimiter{\abs}{\lvert}{\rvert}

\mathtoolsset{showonlyrefs=true}

\newcommand{\cost}{F}

\makeatletter
\newcommand\notsotiny{\@setfontsize\notsotiny\@vipt\@viipt}
\makeatother

\makeatletter
\renewcommand\section{%
  \@startsection{section}{1}{\z@}%
    {1.7ex plus .2ex minus .2ex}%
    {1.3ex plus .1ex minus .1ex}%
    {\normalfont\large\bfseries}%
}
\renewcommand\subsection{%
  \@startsection{subsection}{2}{\z@}%
    {1.0ex plus .2ex minus .2ex}%
    {0.7ex plus .1ex minus .1ex}%
    {\normalfont\normalsize\bfseries}%
}
\makeatother

\makeatletter
\def\thm@space@setup{%
  \thm@preskip=2pt plus 2pt minus 0pt %
  \thm@postskip=2pt plus 2pt minus 0pt %
}
\makeatother

\usepackage{etoolbox}
\makeatletter
\newcommand{\tightdisplays}{%
  \setlength{\abovedisplayskip}{3pt plus 2pt minus 0pt}%
  \setlength{\belowdisplayskip}{3pt plus 2pt minus 0pt}%
  \setlength{\abovedisplayshortskip}{3pt plus 2pt minus 0pt}%
  \setlength{\belowdisplayshortskip}{3pt plus 2pt minus 0pt}%
  \setlength{\jot}{2pt}%
}
\appto\normalsize{\tightdisplays}
\appto\small{\tightdisplays}
\appto\footnotesize{\tightdisplays}
\appto\scriptsize{\tightdisplays}
\makeatother

\usepackage[font=small,labelfont=bf]{caption}
\usepackage{enumitem}
\setlist{nosep,leftmargin=*,itemsep=0.25ex}

\makeatletter
\let\oldthebibliography\thebibliography
\let\endoldthebibliography\endthebibliography
\renewenvironment{thebibliography}[1]{%
  \oldthebibliography{#1}%
  \setlength{\parskip}{0pt}%
  \setlength{\itemsep}{0pt plus 0.3ex}%
}{\endoldthebibliography}
\makeatother

\title{MINIMIZATION OF NONSMOOTH WEAKLY CONVEX FUNCTION OVER PROX-REGULAR SET FOR ROBUST LOW-RANK MATRIX RECOVERY}

\name{Keita Kume and Isao Yamada\thanks{This work was supported by JSPS Grants-in-Aid (24K23885).}}
\address{Dept. of Information and Communications Engineering, Institute of Science Tokyo, JAPAN
\\Email:\{kume,isao\}@sp.ict.e.titech.ac.jp}
\begin{document}
\ninept
\renewcommand{\baselinestretch}{0.95}\selectfont
\maketitle
\begin{abstract}
  We propose a {\em prox-regular}-type low-rank constrained nonconvex nonsmooth optimization model for Robust Low-Rank Matrix Recovery (RLRMR), i.e., estimate problem of low-rank matrix from an observed signal corrupted by outliers.
  For RLRMR, the $\ell_{1}$-norm has been utilized as a convex loss to detect outliers as well as to keep tractability of optimization models.
  Nevertheless, the $\ell_{1}$-norm is not necessarily an ideal robust loss because
  the $\ell_{1}$-norm tends to overpenalize entries corrupted by outliers of large magnitude.
  In contrast, the proposed model can employ a weakly convex function as a more robust loss, against outliers, than the
  $\ell_{1}$-norm.
  For the proposed model, we present (i) a projected variable smoothing-type algorithm applicable for the minimization of a nonsmooth weakly convex function over a prox-regular set, and (ii) a convergence analysis of the proposed algorithm in terms of stationary point.
  Numerical experiments demonstrate the effectiveness of the proposed model compared with the existing models that employ the $\ell_{1}$-norm.
\end{abstract}
\section{Introduction}
\label{sec:intro}
We consider the task of estimating a low-rank matrix
$\bm{X}^{\star} \in \mathbb{R}^{n_{1}\times n_{2}}$
from an observation
$\bm{y} \coloneqq \mathcal{A}(\bm{X}^{\star})+\bm{\epsilon}\in \mathbb{R}^{m}$
with a linear operator
$\mathcal{A}:\mathbb{R}^{n_{1}\times n_{2}} \to \mathbb{R}^{m}$
and a Gaussian noise
$\bm{\epsilon} \in \mathbb{R}^{m}$~\cite{Li-Zhu-Man-Vidal20,Charisopoulos-Chen-Davis-Diaz-Ding-Drusvyatskiy21,Xu-Li-Zheng24}.
This task arises in, e.g., 
phase retrieval~\cite{Shechtman-Eldar-Cohen-Chapman-Miao-Segev15,Zhen-Ma-Xue24,Yazawa-Kume-Yamada25},
robust principal component analysis~\cite{Candes-Li-Ma-Wright11,Chen-Ge-Jiang-Li25}, and low-rank matrix completion~\cite{Candes-Recht12,Cambier-Absil16,Wang-So-Zoubir23,Yao-Dai25}.
In such applications, the observation
$\bm{y}$
may be corrupted due to outliers caused by impulse/sparse noise unavoidably in the observation process.
In order to mitigate the influence of outliers in the estimation process,
Robust Low-Rank Matrix Recovery (RLRMR) has attracted a great attention in the fields of signal processing and machine learning~\cite{Candes-Li-Ma-Wright11,Chen-Ge-Jiang-Li25,Zhen-Ma-Xue24,Cambier-Absil16,Li-Zhu-Man-Vidal20,Charisopoulos-Chen-Davis-Diaz-Ding-Drusvyatskiy21,Wang-So-Zoubir23,Yazawa-Kume-Yamada25,Xu-Li-Zheng24}.
RLRMR is formulated as follows.

\begin{problem}[RLRMR] \label{problem:RLRMR}
Let an observation
$\bm{y} \in \mathbb{R}^{m}$
satisfy
\begin{equation}
  (i=1,2,3,\ldots,m) \quad
  [\bm{y}]_{i} \coloneqq
  \begin{cases}
    \mathcal{A}_{i}(\bm{X}^{\star}) + [\bm{\epsilon}]_{i}, & \mathrm{if}\ i \in \mathcal{I}_{\rm in};   \\
    \xi_{i},                                               & \mathrm{if}\ i \in \mathcal{I}_{\rm out},
  \end{cases} \hspace{-1em}\label{eq:observation}
\end{equation}
where
$\bm{\epsilon} \in \mathbb{R}^{m}$
is a Gaussian noise,
$\mathcal{I}_{\rm in}, \mathcal{I}_{\rm out}\subset \{1,2,\ldots,m\}$
denote unknown disjoint index sets of inliers and outliers such that
$\mathcal{I}_{\rm in} \cup \mathcal{I}_{\rm out} = \{1,2,\ldots,m\}$
and
$\mathcal{I}_{\rm in} \cap \mathcal{I}_{\rm out} = \emptyset$,
each
$\xi_{i} \in \mathbb{R}\ (i\in\mathcal{I}_{\rm out})$
denotes outlier,
and each
$\mathcal{A}_{i}:\mathbb{R}^{n_{1}\times n_{2}}\to\mathbb{R}\ (i=1,2,\ldots,m)$
is a known linear operator.
For convenience, let
$\mathcal{A}:\mathbb{R}^{n_{1}\times n_{2}}\to\mathbb{R}^{m}:\bm{X}\mapsto [\mathcal{A}_{1}(\bm{X})\ \mathcal{A}_{2}(\bm{X})\cdots \mathcal{A}_{m}(\bm{X})]^{\TT}$.
Then,
\begin{equation}
  \mathrm{recover}\ 
  \bm{X}^{\star} \in \mathcal{L}_{r}\subset \mathbb{R}^{n_{1}\times n_{2}}\ 
  \mathrm{from}\ 
  \bm{y}\in\mathbb{R}^{m} \ \mathrm{in~\eqref{eq:observation}},
  \label{eq:RLRMR}
\end{equation}
where
{
\thickmuskip=0.2\thickmuskip
\medmuskip=0.2\medmuskip
\thinmuskip=0.2\thinmuskip
\arraycolsep=0.2\arraycolsep
$\mathcal{L}_{r}  \coloneqq \{\bm{X} \in \mathbb{R}^{n_{1}\times n_{2}} \mid \mathrm{Rank}(\bm{X}) \leq r\}$
with
$r \ll \min\{n_{1},n_{2}\}$.
}
\end{problem}
For Problem~\ref{problem:RLRMR}, the convex
$\ell_{1}$-norm
$\norm{\cdot}_{1}$
has been utilized as standard loss to detect outliers~\cite{Cambier-Absil16,Li-Zhu-Man-Vidal20,Charisopoulos-Chen-Davis-Diaz-Ding-Drusvyatskiy21,Xu-Li-Zheng24}.
For example,
an optimization model
with an expression of
$\mathcal{L}_{r} = \{\bm{U}\bm{V}^{\TT} \in \mathbb{R}^{n_{1}\times n_{2}} \mid \bm{U} \in \mathbb{R}^{n_{1}\times r}, \bm{V} \in \mathbb{R}^{n_{2}\times r}\}$
and
$\lambda \geq 0$:
\begin{equation}
  \thickmuskip=0.2\thickmuskip
  \medmuskip=0.2\medmuskip
  \thinmuskip=0.2\thinmuskip
  \arraycolsep=0.2\arraycolsep
  \mathop{\mathrm{minimize}}\limits_{\bm{U}\in \mathbb{R}^{n_{1}\times r}, \bm{V}\in\mathbb{R}^{n_{2} \times r}} \norm{\bm{y} - \mathcal{A}(\bm{U}\bm{V}^{\TT})}_{1} + \lambda \norm{\bm{U}^{\TT}\bm{U} - \bm{V}^{\TT}\bm{V}}_{F} \label{eq:RLRMR_decomp}
\end{equation}
has been proposed~\cite{Li-Zhu-Man-Vidal20,Charisopoulos-Chen-Davis-Diaz-Ding-Drusvyatskiy21}
for Problem~\ref{problem:RLRMR}, where
$\norm{\cdot}_{F}$
denotes the Frobenius norm, and
the second term reduces the scaling ambiguities of
$\bm{U} \in \mathbb{R}^{n_{1}\times r}$
and
$\bm{V} \in \mathbb{R}^{n_{2}\times r}$
as
$\bm{U}\bm{V}^{\TT} = (\bm{U}\bm{Q})(\bm{V}\bm{Q}^{-\TT})^{\TT}$
with invertible matrices
$\bm{Q} \in \mathbb{R}^{r\times r}$
(see also~\eqref{eq:NRLAD} for another optimization model).
The cost function in~\eqref{eq:RLRMR_decomp} enjoys key properties, namely, {\em sharpness} and {\em weak convexity}~\cite[Prop. 5, and 6]{Li-Zhu-Man-Vidal20}, for applying a {\em subgradient method}~\cite{Davis-Drusvyatskiy-MacPhee18} to the model~\eqref{eq:RLRMR_decomp}.
Theoretical results toward an exact recovery, and numerical experiments in~\cite{Li-Zhu-Man-Vidal20} suggest the effectiveness of the model~\eqref{eq:RLRMR_decomp} for Problem~\ref{problem:RLRMR}.

Nevertheless, the
$\ell_{1}$-norm is not necessarily an ideal loss for robust signal recovery.
To examine this, consider a case where some outliers
$\xi_{\widehat{i}}$
for
$\widehat{i}\in \mathcal{I}_{\rm out}$
deviate excessively from true measurements
$\mathcal{A}_{\widehat{i}}(\bm{X}^{\star})$.
In this case,
$\abs{[\bm{y}]_{\widehat{i}}-\mathcal{A}_{\widehat{i}}(\bm{X}^{\star})} = \abs{\xi_{\widehat{i}}-\mathcal{A}_{\widehat{i}}(\bm{X}^{\star})}$
becomes large, and
$\norm{\bm{y}-\mathcal{A}(\bm{X}^{\star})}_{1}=\sum_{i=1}^{m}\abs{[\bm{y}]_{i} - \mathcal{A}_{i}(\bm{X}^{\star})}$
is dominated by these few terms.
Consequently,
(i)
a minimizer of
$\norm{\bm{y}-\mathcal{A}(\cdot)}_{1}$
may deviate significantly from
$\bm{X}^{\star}$
to reduce the large residuals
$[\bm{y}]_{\widehat{i}}-\mathcal{A}_{\widehat{i}}(\cdot)$;
(ii)
it is desirable to use a robust loss that saturates
such large residuals
$[\bm{y}]_{\widehat{i}}-\mathcal{A}_{\widehat{i}}(\bm{X}^{\star})$.
  
  {\em Weakly convex functions}, e.g.,
the smoothly clipped absolute deviation (SCAD)~\cite{Fan-Li01} (see also~\eqref{eq:SCAD}) and the Minimax Concave Penalty (MCP)~\cite{Zhang10}, are utilized as promising robust loss functions, where a function
$g:\mathbb{R}^{m}\to\mathbb{R}$
is said to be {\em $\eta$-weakly convex} if
$g(\cdot)+\frac{\eta}{2}\norm{\cdot}^{2}$
is convex with
$\eta > 0$.
Indeed, it has been reported that the use of weakly convex loss
improves estimation performance in the literature of robust signal recovery~\cite{Yang-Shen-Ma-Chen-Gu-So19,Suzuki-Yukawa21}.
Hence, the following modification of~\eqref{eq:RLRMR_decomp} with a weakly convex loss
$\ell:\mathbb{R}\to\mathbb{R}$:
\begin{equation}
  \hspace{-1.5em}
  \thickmuskip=0.1\thickmuskip
  \medmuskip=0.1\medmuskip
  \thinmuskip=0.1\thinmuskip
  \arraycolsep=0.1\arraycolsep
  \mathop{\mathrm{minimize}}\limits_{\bm{U}\in \mathbb{R}^{n_{1}\times r}, \bm{V}\in\mathbb{R}^{n_{2} \times r}} \sum_{i=1}^{m}\ell\left([\bm{y}]_{i} - \mathcal{A}_{i}(\bm{U}\bm{V}^{\TT})\right) + \lambda \norm{\bm{U}^{\TT}\bm{U} - \bm{V}^{\TT}\bm{V}}_{F} \hspace{-0.3em} \label{eq:RLRMR_decomp_weakly}
\end{equation}
is expected to improve the estimation performance for Problem~\ref{problem:RLRMR}.
However, for the model~\eqref{eq:RLRMR_decomp_weakly}, any reliable optimization algorithm has not been reported so far mainly because the cost function in~\eqref{eq:RLRMR_decomp_weakly} with a weakly convex
$\ell$
is not necessarily weakly convex.

In this paper, for Problem~\ref{problem:RLRMR}, we propose an optimization model with a weakly convex loss
$\ell:\mathbb{R}\to \mathbb{R}$:
\begin{equation}
  \mathop{\mathrm{minimize}} \ \sum_{i=1}^{m}\ell\left([\bm{y}]_{i} - \mathcal{A}_{i}(\bm{X})\right) \ \mathrm{subject\ to\ }
  \bm{X}\in \mathfrak{L}_{r,\sigma},
  \label{eq:proposed_model}
\end{equation}
where
$\mathfrak{L}_{r,\sigma} \subset \mathfrak{L}_{r}\subset \mathbb{R}^{n_{1}\times n_{2}}$
with
$\sigma > 0$
is a closed nonconvex set:
\begin{equation}
  \hspace{-1.2em}
  \mathfrak{L}_{r,\sigma}\coloneqq \left\{\bm{X} \in \mathbb{R}^{n_{1}\times n_{2}} \mid \substack{0 < \mathrm{Rank}(\bm{X}) \leq r,\\ \sigma \leq \sigma_{j}(\bm{X})\ \mathrm{or}\ \sigma_{j}(\bm{X})=0 \ (j=1,2,\ldots, r)}\right\},\hspace{-0.4em} \label{eq:low_rank}
\end{equation}
$\sigma_{j}(\cdot)$
stands for the $j$th largest singular value of a given matrix,
and
$\sigma > 0$
works as a lower threshold for nonzero singular values.
The set
$\mathfrak{L}_{r,\sigma}$
with a small
$\sigma > 0$
is a reasonable approximation of
$\mathfrak{L}_{r}$
because
(i)
$\mathfrak{L}_{r, \sigma}\cup \{\bm{0}\}$
converges to
$\mathfrak{L}_{r}$
as
$\sigma \searrow 0$
in the sense of Painlev\'e-Kuratowski~\cite[Sect. 4.B]{Rockafellar-Wets98};
and
(ii)
$\mathfrak{L}_{r,\sigma}$
with
$\sigma > 0$
is {\em prox-regular}~\cite[Thm. 5]{Balashov-Kamalov21} (see Problem~\ref{problem:origin}~\ref{enum:problem:origin:constraint}), which serves as key properties (see, e.g.,~\cite{Clarke-Stern-Wolenski95,Poliquin-Rockafellar-Thibault00,Rockafellar-Wets98})
regarding the {\em metric projection} onto
$\mathfrak{L}_{r,\sigma}$,
where
$\mathfrak{L}_{r}$
is not prox-regular~\cite{Levin-Kileel-Boumal24}.
Moreover, the cost function in~\eqref{eq:proposed_model} with a weakly convex
$\ell$
remains weakly convex
as a composition of
$\ell$
with an affine operator.
Hence, the model~\eqref{eq:proposed_model} seems to be more tractable than the model~\eqref{eq:RLRMR_decomp_weakly}.
To the best of the authors' knowledge, this is the first work to employ
$\mathcal{L}_{r,\sigma}$
as the constraint set in optimization models for Problem~\ref{problem:RLRMR}.

We also present an iterative algorithm for the model~\eqref{eq:proposed_model} via the following optimization problem over a prox-regular set.
Indeed, the model~\eqref{eq:proposed_model} is reformulated into Problem~\ref{problem:origin} (see Remark~\ref{remark:formulation}).
\begin{problem}\label{problem:origin}
Let
$\mathcal{X}$
and
$\mathcal{Z}$
be Euclidean spaces.
Then,
\begin{equation}
  \mathrm{find} \ \bm{x}^{\star} \in \argmin_{\bm{x} \in C} \cost(\bm{x}) (\neq \emptyset), \label{eq:problem}
\end{equation}
where
$\cost\coloneqq g\circ\mathfrak{S}$,
$\mathfrak{S}:\mathcal{X}\to\mathcal{Z}$,
$g:\mathcal{Z}\to\mathbb{R}$
and
$C\subset\mathcal{X}$
satisfy 
\begin{enumerate}[label=(\roman*)]
  \item
        $\mathfrak{S}: \mathcal{X} \to \mathcal{Z}$
        is a differentiable mapping such that
        $\mathfrak{S}$
        and
        its G\^{a}teaux derivative
        $\mathrm{D}\mathfrak{S}$
        are Lipschitz continuous;
  \item
        $g:\mathcal{Z} \to \mathbb{R}$
        is (a) $L_{g}$-Lipschitz continuous with
        $L_{g} > 0$
        (possibly nonsmooth),
        (b)~$\eta$-weakly convex with
        $\eta > 0$,
        i.e.,
        $g+\frac{\eta}{2}\norm{\cdot}^{2}$
        is convex,
        and (c) {\em prox-friendly}, i.e.,
        $\prox{\mu g}\ (\mu \in (0,\eta^{-1}))$
        (see~\eqref{eq:prox})
        is available as a computable tool;
  \item \label{enum:problem:origin:constraint}
        $C\subset\mathcal{X}$
        is a nonempty closed {\em prox-regular set}\footnote{
          Other than
          $\mathfrak{L}_{r,\sigma}$
          in~\eqref{eq:low_rank},
          prox-regular sets include, e.g., closed convex sets,
          $C^{2}$ embedded submanifolds in
          $\mathcal{X}$ (see, e.g.,~\cite{Federer59}, and~\cite[Lemma 2.1]{Lewis-Malick08}), and
          $\{\bm{x} \in \mathbb{R}^{n} \mid \norm{\bm{x}}=1, [\bm{x}]_{i} \geq 0\ \forall i\}$ (see~\cite[Lemma 3 (ii)]{Hu-Deng-Wu-Li24}).
        }~\cite[Thm. 1.3]{Poliquin-Rockafellar-Thibault00}, i.e.,
        the {\em metric projection}
        $P_{C}:\mathcal{X} \rightrightarrows C:\widebar{\bm{x}}\mapsto \mathop{\mathrm{argmin}}_{\bm{x}\in C}\norm{\widebar{\bm{x}}-\bm{x}}$
        onto
        $C$
        is single-valued on some open superset of
        $C$.
        Moreover, we assume that at least
        one point
        $\bm{x} \in P_{C}(\widebar{\bm{x}})$
        can be computed for every
        $\widebar{\bm{x}} \in \mathcal{X}$.
\end{enumerate}
\end{problem}

\begin{remark}[Reformulation of the model~\eqref{eq:proposed_model} into Problem~\ref{problem:origin}] \label{remark:formulation}
  The proposed model~\eqref{eq:proposed_model} is a special instance of Problem~\ref{problem:origin} by setting
  $\mathcal{X}\coloneqq\mathbb{R}^{n_{1}\times n_{2}}$,
  $\mathcal{Z}\coloneqq \mathbb{R}^{m}$,
  $\mathfrak{S}(\bm{X})\coloneqq \bm{y}-\mathcal{A}(\bm{X})$,
  $g:\mathbb{R}^{m}\to\mathbb{R}:\bm{z} \mapsto \sum_{i=1}^{m}\ell([\bm{z}]_{i})$,
  and
  $C\coloneqq \mathfrak{L}_{r,\sigma}$
  in~\eqref{eq:low_rank} with
  $\sigma > 0$,
  where
  $\ell$
  is Lipschitz continuous, weakly convex and prox-friendly
  (e.g., the $\ell_{1}$-norm,
  SCAD~\cite{Fan-Li01}, and MCP~\cite{Zhang10}. See also~\cite[Sect. 2.2]{Kume-Yamada24A} for other examples),
  and
  $P_{\mathcal{L}_{r,\sigma}}$
  can be computed by using a singular value decomposition (SVD) of a given matrix~\cite[Cor. 2]{Balashov-Kamalov21}.
\end{remark}

In order to find a {\em stationary point} of Problem~\ref{problem:origin},
we extend variable smoothing-type algorithms~\cite{Liu-Xia24,Lopes-Peres-Bilches25,Kume-Yamada25A,Kume-Yamada25B} so that even a nonconvex prox-regular set, e.g.,
$\mathfrak{L}_{r,\sigma}$,
can be used as
$C$,
whereas only a closed convex set
$C$
can be used in~\cite{Liu-Xia24,Lopes-Peres-Bilches25,Kume-Yamada25A,Kume-Yamada25B}.
The proposed algorithm, named {\em projected variable smoothing}, is designed as a projected gradient method with a smoothed surrogate function
$\moreau{g}{\mu}\circ\mathfrak{S}$
of
$g\circ\mathfrak{S}$,
where
$\moreau{g}{\mu}:\mathcal{Z}\to\mathbb{R}$
with
$\mu \in (0,\eta^{-1})$
is the {\em Moreau envelope} of
$g$ (see~\eqref{eq:moreau}).
By exploiting its notable properties, e.g.,
$\lim_{\mu\to 0}\moreau{g}{\mu}(\bm{z})= g(\bm{z})\ (\bm{z}\in \mathcal{Z})$,
the proposed algorithm updates
$(\bm{x}_{n})_{n=1}^{\infty} \subset C$
as
$\bm{x}_{n+1} \in P_{C}(\bm{x}_{n}-\gamma_{n}\nabla (\moreau{g}{\mu_{n}}\circ\mathfrak{S})(\bm{x}_{n}))$,
where
$\gamma_{n} > 0$
and
$\mu_{n} (\in (0,\eta^{-1})) \searrow 0$
are chosen strategically.
We also present an asymptotic convergence analysis, in Theorem~\ref{theorem:convergence_extension}, of the proposed algorithm in terms of a {\em stationary point}
(see just after~\eqref{eq:Fermat}).

Numerical experiments demonstrate the effectiveness of proposed model~\eqref{eq:proposed_model}  for Problem~\ref{problem:RLRMR} solved by the proposed algorithm.

  {\bf Related work on Problem~\ref{problem:origin}}:
The theory of prox-regular sets dates back to 1959~\cite{Federer59}, and has been studied in, e.g.,~\cite{Clarke-Stern-Wolenski95,Poliquin-Rockafellar-Thibault00}.
However,
to the best of our knowledge, only a few recent papers~\cite{Balashov20,Davis-Drusvyatskiy-Shi25} propose iterative algorithms for Problem~\ref{problem:origin}
under the assumption that
$C$
is {\em proximally smooth}\footnote{
  $C$
  is said to be {\em proximally smooth} if
  $P_{C}$
  is single-valued over
  $C + B(\bm{0}; \delta) \coloneqq \{\bm{x}+\bm{u} \in \mathcal{X}\mid \bm{x} \in C, \bm{u}\in\mathcal{X}, \norm{\bm{u}} < \delta\}$
  with some
  $\delta > 0$.
}, which is stronger than prox-regularity.
Specifically, \cite{Balashov20} extends
a classical projected subgradient method to the setting of Problem~\ref{problem:origin}, and presents its linear-rate convergence under an additional {\em error-bound condition} on
$\cost$.
The paper~\cite{Davis-Drusvyatskiy-Shi25} extends classical algorithms, e.g., a projected subgradient method and a proximal point method, to the setting of Problem~\ref{problem:origin}, and presents iteration-complexity bounds to reach a predetermined tolerance.
However, any asymptotic convergence analysis of algorithms in~\cite{Davis-Drusvyatskiy-Shi25} has not been reported yet.
In contrast, the proposed algorithm has an asymptotic convergence guarantee (see Thm.~\ref{theorem:convergence_extension}) for Problem~\ref{problem:origin} without error-bound condition required in~\cite{Balashov20}.

{
\thickmuskip=0.0\thickmuskip
\medmuskip=0.0\medmuskip
\thinmuskip=0.0\thinmuskip
\arraycolsep=0.0\arraycolsep
{\bf Notation}:
$\mathbb{N}$,
$\mathbb{R}$,
and
$\mathbb{R}_{++}$,
denote respectively the sets of all positive integers, all real numbers, and all positive real numbers.
$[\bm{v}]_{i} \in \mathbb{R}$
denotes
$i$th entry of
$\bm{v} \in \mathbb{R}^{n}$,
For
$\widebar{\bm{x}} \in \mathcal{X}$
and a nonempty closed
set
$E \subset \mathcal{X}$,
$\|\widebar{\bm{x}}\|\coloneqq \sqrt{\inprod{\widebar{\bm{x}}}{\widebar{\bm{x}}}}$
denotes the Euclidean norm with the standard inner product
$\inprod{\cdot}{\cdot}$,
and
$\dist(\widebar{\bm{x}},E)\coloneqq \min\{\norm{\bm{v}-\widebar{\bm{x}}}\mid \bm{v} \in E\}$
denotes the distance function.
For a differentiable mapping
$\mathcal{F}:\mathcal{X}\to \mathcal{Z}$,
its G\^{a}teaux derivative at
$\widebar{\bm{x}}\in \mathcal{X}$
is the linear operator
$\mathrm{D}\mathcal{F}(\widebar{\bm{x}}):\mathcal{X} \to \mathcal{Z}:\bm{v}\mapsto \lim_{\mathbb{R}\setminus\{0\} \ni t\to 0}\frac{\mathcal{F}(\widebar{\bm{x}}+t\bm{v}) -\mathcal{F}(\widebar{\bm{x}})}{t}$.
A mapping
$\mathcal{F}:\mathcal{X}\to\mathcal{Z}$
is said to be {\em Lipschitz continuous} with a Lipschitz constant
$L_{\mathcal{F}} > 0$
if
$\norm{\mathcal{F}(\bm{x}_{1}) - \mathcal{F}(\bm{x}_{2})} \leq L_{\mathcal{F}}\norm{\bm{x}_{1}-\bm{x}_{2}}\ (\bm{x}_{1},\bm{x}_{2} \in \mathcal{X})$.
For a differentiable function
$J:\mathcal{X} \to \mathbb{R}$,
$\nabla J(\widebar{\bm{x}}) \in \mathcal{X}$
is the gradient of
$J$
at
$\widebar{\bm{x}} \in \mathcal{X}$
if
$\mathrm{D}J(\widebar{\bm{x}})[\bm{v}] = \inprod{\nabla J(\widebar{\bm{x}})}{\bm{v}}\ (\bm{v} \in \mathcal{X})$.
A function
$J:\mathcal{X} \to \exR$
is said to be {\em proper} if
$\dom{J}\coloneqq\{\bm{x}\in \mathcal{X} \mid J(\bm{x}) <\infty\} \neq \emptyset$.
}

\section{Preliminary on Nonsmooth Analysis}
We review necessary notions and tools in nonsmooth analysis in~\cite{Rockafellar-Wets98} (see also a comprehensive review paper~\cite{Li-So-Ma20}).
\begin{definition}[Subdifferential~{\cite[Def. 8.3]{Rockafellar-Wets98}}]\label{definition:subdifferential}
  A vector
  $\bm{v} \in \mathcal{X}$
  is said to be a {\em Fr\'echet (or regular) subgradient} of a proper function
  $J:\mathcal{X} \to \exR$
  at
  $\widebar{\bm{x}} \in \dom{J}$,
  denoted by
  $\bm{v} \in \Fsubdiff J(\widebar{\bm{x}})$~\cite[Def. 8.3]{Rockafellar-Wets98},
  if
  $
    \sup\limits_{\epsilon>0} \left(\inf\limits_{0<\norm{\bm{x}-\widebar{\bm{x}}}< \epsilon}\frac{J(\bm{x})-J(\widebar{\bm{x}}) - \inprod{\bm{v}}{\bm{x}-\widebar{\bm{x}}}}{\norm{\bm{x}-\widebar{\bm{x}}}}\right) \geq 0$
  holds,
  where
  $\Fsubdiff J$
  is called the {\em Fr\'echet subdifferential} of
  $J$,
  and
  $\Fsubdiff J(\widebar{\bm{x}})$
  at
  $\widebar{\bm{x}} \notin \dom{J}$
  is understood as
  $\emptyset$.
  If
  $J$
  is convex, then
  $\Fsubdiff J$
  coincides with the convex subdifferential of
  $J$~\cite[Prop. 8.12]{Rockafellar-Wets98}.
\end{definition}

Fermat's rule~\cite[Thm. 10.1]{Rockafellar-Wets98} serves as a necessary condition:
\begin{equation}
  \Fsubdiff (\cost+\iota_{C})(\bm{x}^{\star})  \ni \bm{0}, \label{eq:Fermat}
\end{equation}
for local optimality of Problem~\ref{problem:origin},
where
the {\em indicator function}
$\iota_{C}:\mathcal{X}\to\exR$
is defined as
$\iota_{C}(\widebar{\bm{x}}) \coloneqq 0$
if
$\widebar{\bm{x}} \in C$;
$\iota_{C}(\widebar{\bm{x}}) \coloneqq +\infty$
if
$\widebar{\bm{x}} \notin C$.
A point
$\bm{x}^{\star}$
enjoying~\eqref{eq:Fermat} is called a {\em stationary point} of Problem~\ref{problem:origin}, and finding a stationary point is a reasonable goal for nonconvex optimization~\cite{Rockafellar-Wets98,Davis-Drusvyatskiy-MacPhee18,Li-So-Ma20,Li-Zhu-Man-Vidal20,Bohm-Wright21,Liu-Xia24,Kume-Yamada24,Kume-Yamada24A,Lopes-Peres-Bilches25,Kume-Yamada25A,Kume-Yamada25B}.
In this paper, we aim to find a stationary point of Problem~\ref{problem:origin}.

{\em The proximity operator} and {\em the Moreau envelope} have been used as computational tools for nonsmooth optimization~\cite{Yamada-Yukawa-Yamagishi11,Bohm-Wright21,Liu-Xia24,Kume-Yamada24,Kume-Yamada24A,Lopes-Peres-Bilches25,Kume-Yamada25A,Kume-Yamada25B,Yazawa-Kume-Yamada25}.
For a Euclidean space
$\mathcal{H}$,
the proximity operator of a proper function
$J:\mathcal{H}\to\exR$
with index
$\mu > 0$
is defined by
\begin{equation}
  \hspace{-0.5em}
  \prox{\mu J}:\mathcal{H} \rightrightarrows \mathcal{H}:
  \widebar{\bm{u}} \mapsto
  \argmin_{\bm{u} \in \mathcal{H}} \left(J(\bm{u}) + \frac{1}{2\mu}\|\bm{u}-\widebar{\bm{u}}\|^{2}\right). \label{eq:prox}
\end{equation}
By letting
$\mathcal{H}\coloneqq \mathcal{X}$
and
$J\coloneqq \iota_{C}$
in Problem~\ref{problem:origin}, we have the expression
$P_{C} = \prox{\mu\iota_{C}}\ (\mu > 0)$,
and
$P_{C}(\widebar{\bm{x}}) \subset C$
is a nonempty and compact set for every
$\widebar{\bm{x}} \in \mathcal{X}$~\cite[Thm. 1.25]{Rockafellar-Wets98}.
On the other hand, by letting
$\mathcal{H}\coloneqq \mathcal{Z}$
and
$J\coloneqq g$
in Problem~\ref{problem:origin},
$\prox{\mu g}(\widebar{\bm{z}}) \in \mathcal{Z}$
with
$\mu \in (0,\eta^{-1})$
is single-valued for every
$\widebar{\bm{z}} \in \mathcal{Z}$
due to the strong convexity of
$g + (2\mu)^{-1}\|\cdot-\widebar{\bm{z}}\|^{2}$.
By using
$\prox{\mu g}$,
the Moreau envelope
$\moreau{g}{\mu}:\mathcal{Z}\to\mathbb{R}\ (\mu \in (0,\eta^{-1}))$
of
$g$
is defined by
\begin{equation}
  \moreau{g}{\mu}:\mathcal{Z} \to \mathbb{R}:\widebar{\bm{z}} \mapsto
  g(\prox{\mu g}(\widebar{\bm{z}})) + \frac{1}{2\mu}\norm{\prox{\mu g}(\widebar{\bm{z}}) - \widebar{\bm{z}}}^{2}.
  \label{eq:moreau}
\end{equation}
$\moreau{g}{\mu}$
has notable properties as an approximation of the nonsmooth function
$g$:
(i)
$\lim_{\mu\to 0}\moreau{g}{\mu}(\widebar{\bm{z}}) = g(\widebar{\bm{z}})\ (\widebar{\bm{z}}\in\mathcal{Z})$;
(ii) continuous differentiability of
$\moreau{g}{\mu}$
with
$\nabla\moreau{g}{\mu}(\widebar{\bm{z}})=\mu^{-1}(\widebar{\bm{z}} - \prox{\mu g}(\widebar{\bm{z}}))$;
(iii)
Lipschitz continuity of
$\nabla\moreau{g}{\mu}$
(see, e.g,~\cite[Cor. 3.4]{Hoheisel-Laborde-Oberman20}).
Both
$\moreau{g}{\mu}$
and
$\nabla \moreau{g}{\mu}$
are computable if
$g$
is prox-friendly (see Rem.~\ref{remark:formulation} for such
$g$).
Moreover, 
$\nabla (\moreau{g}{\mu}\circ\mathfrak{S})$
enjoys the Lipschitz continuity.
\begin{fact}[Lipschitzian of $\nabla (\moreau{g}{\mu}\circ\mathfrak{S})${~\cite[Prop. 4.2(a)]{Kume-Yamada24A}}]
  \label{fact:Lipschitz}
  Consider
  $g\circ\mathfrak{S}$
  in Problem~\ref{problem:origin}.
  Then,
  $\nabla (\moreau{g}{\mu}\circ\mathfrak{S})$
  is Lipschitz continuous with a Lipschitz constant
  $L_{\nabla (\moreau{g}{\mu}\circ\mathfrak{S})} > 0$,
  where there exist
  $\varpi_{1},\varpi_{2} \in \mathbb{R}_{++}$
  such that
  $L_{\nabla (\moreau{g}{\mu}\circ\mathfrak{S})} = \varpi_{1} + \varpi_{2}\mu^{-1}\ (\forall \mu \in (0,(2\eta)^{-1}])$.
\end{fact}

\section{Projected variable smoothing algorithm}
This section presents a key idea, based on a stationarity measure, for finding a stationary point of Problem~\ref{problem:origin} (see~\eqref{eq:Fermat}), and then an iterative algorithm  for Problem~\ref{problem:origin} with a convergence analysis.
Due to space limitation, all proofs of the theoretical results in this paper
are deferred to an extended manuscript in preparation.

Inspired by~\cite{Liu-Xia24,Kume-Yamada25A,Kume-Yamada25B},
we employ
a {\em gradient mapping-type stationarity measure}
$\mathcal{M}_{\gamma}^{\cost,\iota_{C}}:\mathcal{X}\to \mathbb{R}\ (\gamma\in\mathbb{R}_{++})$
defined
as
\begin{equation}
  \thickmuskip=0.3\thickmuskip
  \medmuskip=0.3\medmuskip
  \thinmuskip=0.3\thinmuskip
  \arraycolsep=0.3\arraycolsep
  (\widebar{\bm{x}} \in \mathcal{X}) \quad
  \mathcal{M}_{\gamma}^{\cost,\iota_{C}}(\widebar{\bm{x}})
  \coloneqq
  \dist\left(\bm{0},\frac{\widebar{\bm{x}} - P_{C}\left(\widebar{\bm{x}}-\gamma \Fsubdiff \cost\left(\widebar{\bm{x}}\right)\right)}{\gamma}\right),
  \label{eq:measure}
\end{equation}
with
$\frac{\widebar{\bm{x}} - P_{C}\left(\widebar{\bm{x}}-\gamma \Fsubdiff \cost\left(\widebar{\bm{x}}\right)\right)}{\gamma}
  = \left\{\frac{\widebar{\bm{x}} - \bm{p}}{\gamma} \mid  \bm{p}\in P_{C}(\widebar{\bm{x}}-\gamma \bm{v}), \ \bm{v} \in \Fsubdiff \cost(\widebar{\bm{x}})\right\}$.
In a special case where
$C$
is closed convex,
some useful properties of
$\mathcal{M}_{\gamma}^{\cost,\iota_{C}}$
are found in~\cite{Liu-Xia24,Kume-Yamada25A,Kume-Yamada25B}.
For example, under the convexity of
$C$,
a stationary point
$\bm{x}^{\star} \in \mathcal{X}$
can be characterized by
$\mathcal{M}_{\gamma}^{\cost,\iota_{C}}(\bm{x}^{\star})=0$
for every
$\gamma \in \mathbb{R}_{++}$~\cite[Fact III.1 (a)]{Kume-Yamada25B}.
Even in the absence of the convexity of
$C$,
we can characterize a stationary point via
$\mathcal{M}_{\gamma}^{\cost,\iota_{C}}$
for some
$\gamma \in \mathbb{R}_{++}$
by Lemma~\ref{lemma:stationarity}.
\begin{lemma}[Stationarity characterization via $\mathcal{M}_{\gamma}^{\cost,\iota_{C}}$] \label{lemma:stationarity}
  $\bm{x}^{\star} \in \mathcal{X}$
  is a stationary point of Problem~\ref{problem:origin}, i.e.,
  $\Fsubdiff(\cost+\iota_{C})(\bm{x}^{\star}) \ni \bm{0}$,
  if and only if
  $\mathcal{M}_{\gamma}^{\cost,\iota_{C}}(\bm{x}^{\star}) = 0$
  holds for some
  $\gamma \in \mathbb{R}_{++}$.
\end{lemma}

By Lemma~\ref{lemma:stationarity}, we can find a stationary point of Problem~\ref{problem:origin} by finding a point
$\bm{x}^{\star} \in \mathcal{X}$
where
$\mathcal{M}_{\gamma}^{\cost,\iota_{C}}(\bm{x}^{\star})=0$
holds for some
$\gamma\in\mathbb{R}_{++}$.
However, it is still challenging to approximate iteratively a point
$\bm{x}^{\star}$
with
$\mathcal{M}_{\gamma}^{\cost,\iota_{C}}(\bm{x}^{\star})=0$
due to the nonsmoothness of
$g$
in
$\cost=g\circ\mathfrak{S}$.
In contrast,
if
$\cost$
were continuously differentiable, then we could exploit powerful arts, e.g., {\em sufficient decrease property},
developed for proximal (or projected) gradient method.

\begin{fact}[Sufficient decrease property, e.g.,~{\cite[Lemma 2]{Bolte-Sabach-Teboulle14}}] \label{fact:Armijo}
  Consider
  $C\subset \mathcal{X}$
  in Problem~\ref{problem:origin}.
  Let
  $J:\mathcal{X}\to\mathbb{R}$
  be continuously differentiable such that
  $\nabla J:\mathcal{X}\to\mathcal{X}$
  is Lipschitz continuous with a Lipschitz constant
  $L_{\nabla J}>0$.
  Let
  $c \in (0,1/2)$.
  Then, we have
    {
      \thickmuskip=0.2\thickmuskip
      \medmuskip=0.2\medmuskip
      \thinmuskip=0.2\thinmuskip
      \arraycolsep=0.2\arraycolsep
      \begin{align}
         & (\forall \gamma \in (0,(1-2c)L_{\nabla J}^{-1}],\ \forall \widebar{\bm{x}}\in C,\ \forall \bm{x} \in P_{C}(\widebar{\bm{x}}-\gamma\nabla J(\widebar{\bm{x}})) \subset C) \\
         & \ J(\bm{x})
        \leq J(\widebar{\bm{x}}) - c\gamma \norm{(\widebar{\bm{x}}-\bm{x})/\gamma}^{2}
        \leq J(\widebar{\bm{x}}) - c\gamma \left(\mathcal{M}_{\gamma}^{J,\iota_{C}}(\widebar{\bm{x}})\right)^{2}. \label{eq:Armijo_general}
      \end{align} }%
\end{fact}

Fact~\ref{fact:Armijo} motivates us to replace the nonsmooth function
$g$
in
$\mathcal{M}_{\gamma}^{\cost,\iota_{C}} = \mathcal{M}_{\gamma}^{g\circ\mathfrak{S},\iota_{C}}$
with its Moreau envelope
$\moreau{g}{\mu}$
in~\eqref{eq:moreau}.
Indeed,
the gradient
$\nabla (\moreau{g}{\mu}\circ\mathfrak{S})$
is Lipschitz continuous by Fact~\ref{fact:Lipschitz},
and thus
$\mathcal{M}_{\gamma}^{\moreau{g}{\mu}\circ\mathfrak{S},\iota_{C}}$
enjoys the sufficient decrease property in~\eqref{eq:Armijo_general} with
$J\coloneqq \moreau{g}{\mu}\circ\mathfrak{S}$.
However, we have a gap between
$\mathcal{M}_{\gamma}^{\moreau{g}{\mu}\circ\mathfrak{S},\iota_{C}}$
and
$\mathcal{M}_{\gamma}^{\cost,\iota_{C}}$
because the condition
$\mathcal{M}_{\gamma}^{\moreau{g}{\mu}\circ\mathfrak{S},\iota_{C}}(\widebar{\bm{x}}) = 0$
for
$\widebar{\bm{x}}\in \mathcal{X}$
does not imply
$\mathcal{M}_{\gamma}^{\cost,\iota_{C}}(\widebar{\bm{x}})=0$.
Theorem~\ref{theorem:asymptotic_approximation} bridges the gap by passing through an asymptotic behavior of
$\mathcal{M}_{\gamma}^{\moreau{g}{\mu}\circ\mathfrak{S},\iota_{C}}$
as
$\mu \searrow 0$.
\begin{theorem}[\footnote{
      In a special case where
      $C$
      is closed convex and
      $\gamma_{n}\coloneqq \widebar{\gamma} > 0\ (n\in\mathbb{N})$,
      a similar inequality to~\eqref{eq:smoothing_asymptotic} is found in our recent paper~\cite[Thm. III.2 (b)]{Kume-Yamada25B}.
      However, an extension of~\cite[Thm. III.2 (b)]{Kume-Yamada25B} to the setting in Thm.~\ref{theorem:asymptotic_approximation}, 
      especially to the variable sequence
      $(\gamma_{n})_{n=1}^{\infty}$,
      is non-trivial.
    }Asymptotic property of  and $\mathcal{M}_{\gamma}^{\moreau{g}{\mu}\circ\mathfrak{S},\phi}$]
  \label{theorem:asymptotic_approximation}
  Consider
  $\cost=g\circ\mathfrak{S}$
  and
  $C$
  in Problem~\ref{problem:origin}.
  Let
  $(\bm{x}_{n})_{n=1}^{\infty} \subset C$
  and
  $(\gamma_{n})_{n=1}^{\infty} \subset \mathbb{R}_{++}$
  converge respectively to some
  $\widebar{\bm{x}} \in C$
  and
  $\widebar{\gamma} \geq 0$.
  Then,
  for
  $\mu_{n}(\in (0,(2\eta)^{-1}]) \searrow 0$
  and
  $\cost_{n} \coloneqq  \moreau{g}{\mu_{n}}\circ\mathfrak{S}\ (n\in\mathbb{N})$,
  \begin{equation}
    \thickmuskip=0.3\thickmuskip
    \medmuskip=0.2\medmuskip
    \thinmuskip=0.2\thinmuskip
    \arraycolsep=0.2\arraycolsep
    \hspace{-1.6em}
    \liminf_{n\to\infty} \mathcal{M}_{\gamma_{n}}^{\cost_{n},\iota_{C}}(\bm{x}_{n})
    \geq
    \begin{cases}
      \mathcal{M}_{\widebar{\gamma}}^{\cost,\iota_{C}}(\widebar{\bm{x}}), & \hspace{-0.8em} \mathrm{if}\ \widebar{\gamma}>0; \\
      \dist\left(\bm{0}, \Fsubdiff (\cost+\iota_{C})(\widebar{\bm{x}})\right), 
                                                                          & \hspace{-0.8em}

      \mathrm{if}\ \widebar{\gamma} = 0
    \end{cases}\hspace{-1em} \label{eq:smoothing_asymptotic}
  \end{equation}
  holds.
  Moreover,
  by combining Lemma~\ref{lemma:stationarity} and~\eqref{eq:Fermat},
  $\widebar{\bm{x}}$
  is a stationary point of Problem~\ref{problem:origin} if
  $\liminf_{n\to\infty} \mathcal{M}_{\gamma_{n}}^{\cost_{n},\iota_{C}}(\bm{x}_{n}) = 0$.
\end{theorem}

By Thm.~\ref{theorem:asymptotic_approximation}, our goal for finding a stationary point of Problem~\ref{problem:origin} is reduced to finding a sequence
$(\bm{x}_{n},\gamma_{n})_{n=1}^{\infty} \subset C \times \mathbb{R}_{++}$
such that
$\liminf_{n\to\infty}\mathcal{M}_{\gamma_{n}}^{\cost_{n} ,\iota_{C}}(\bm{x}_{n}) = 0$
with some
$\mu_{n}(\in (0,(2\eta)^{-1}]) \searrow 0$
and
$\cost_{n} \coloneqq \moreau{g}{\mu_{n}}\circ\mathfrak{S}$.
Based on this observation, we present a projected gradient-type method, called a {\em projected variable smoothing algorithm}, illustrated in Algorithm~\ref{alg:proposed}.

\begin{algorithm}[t]
  {
    \small
    \caption{Projected variable smoothing for Problem~\ref{problem:origin}} \label{alg:proposed}
    {
      \begin{algorithmic}[1]
        \Require
        $\bm{x}_{1}\in C$,
        $c\in(0,2^{-1})$,
        $\rho \in (0,1)$,
        $\widetilde{\gamma} \in \mathbb{R}_{++}$,
        $\alpha \geq 1$
        \For{$n=1,2,3,\ldots$}
        \State
        Set
        $\mu_{n}\leftarrow (2\eta)^{-1}n^{-1/\alpha}$
        and
        $\cost_{n}\coloneqq \moreau{g}{\mu_{n}}\circ\mathfrak{S}$
        \For{$m=0,1,2,\ldots$}\label{line:backtracking_start}
        \State
        $\bm{x} \in P_{C}(\bm{x}_{n} - \rho^{m}\widetilde{\gamma} \nabla \cost_{n}(\bm{x}_{n})) \subset C$ \label{line:projection}
        \If{$\cost_{n}(\bm{x}) \leq \cost_{n}(\bm{x}_{n}) - c\rho^{m}\widetilde{\gamma}\norm{\frac{\bm{x}_{n}-\bm{x}}{\rho^{m}\widetilde{\gamma}}}^{2}$ holds} \label{line:Armijo}
        \State
        Set
        $(\bm{x}_{n+1},\gamma_{n}) \leftarrow (\bm{x},\rho^{m}\widetilde{\gamma})$, and \textbf{break}
        \EndIf
        \EndFor \label{line:backtracking_end}
        \EndFor
      \end{algorithmic}
    }
  }
\end{algorithm}

In Alg.~\ref{alg:proposed}, we update estimates as
$\bm{x}_{n+1} \in P_{C}(\bm{x}_{n}-\gamma_{n}\nabla \cost_{n}(\bm{x}_{n}))$
with a stepsize
$\gamma_{n}$
such that the sufficient decrease condition in~\eqref{eq:Armijo_general} with
$\widebar{\bm{x}} \coloneqq \bm{x}_{n}$
and
$J\coloneqq \moreau{g}{\mu_{n}}\circ\mathfrak{S}$
is achieved by
$(\bm{x},\gamma)\coloneqq (\bm{x}_{n+1},\gamma_{n})$,
where
$\mu_{n}\coloneqq (2\eta)^{-1}n^{-1/\alpha}$
with
$\alpha \geq 1$
(Note: $\alpha=3$ is recommended in the literature of variable smoothing-type algorithms~\cite{Bohm-Wright21,Kume-Yamada24,Kume-Yamada24A,Liu-Xia24,Lopes-Peres-Bilches25,Kume-Yamada25A,Yazawa-Kume-Yamada25,Kume-Yamada25B} that can not be applied to Problem~\ref{problem:origin} due to the nonconvexity of
$C$).
To find such a pair
$(\bm{x}_{n+1},\gamma_{n})$,
we employ a standard {\em backtracking algorithm} in line~\ref{line:backtracking_start}-\ref{line:backtracking_end} of Alg.~\ref{alg:proposed} (see, e.g.,~\cite{Kume-Yamada24A,Kume-Yamada25B}).
Thanks to Fact~\ref{fact:Armijo} with Fact~\ref{fact:Lipschitz},
we can find such
$(\bm{x}_{n+1},\gamma_{n})$
in at most
$\max\left\{1, \left\lceil \log_{\rho}\left(\frac{1-2c}{\widetilde{\gamma}L_{\nabla (\moreau{g}{\mu_{n}}\circ\mathfrak{S})}}\right)\right\rceil\right\}$
backtracking steps in line~\ref{line:backtracking_start}-\ref{line:backtracking_end}, where
Alg.~\ref{alg:proposed} does not require any knowledge on a constant
$L_{\nabla (\moreau{g}{\mu_{n}}\circ\mathfrak{S})}=\varpi_{1}+\varpi_{2}\mu_{n}^{-1}$
in Fact~\ref{fact:Lipschitz}.
Finally, we present a convergence analysis of Alg.~\ref{alg:proposed} in Theorem~\ref{theorem:convergence_extension}.
\begin{theorem}[Convergence analysis of  Alg.~\ref{alg:proposed}]\label{theorem:convergence_extension}
  Consider Problem~\ref{problem:origin}.
  Choose
  $\bm{x}_{1} \in C$,
  $c \in (0,2^{-1})$,
  $\rho \in (0,1)$,
  $\widetilde{\gamma} \in \mathbb{R}_{++}$,
  and
  $\alpha \geq 1$.
  For
  $(\bm{x}_{n})_{n=1}^{\infty} \subset C$
  generated by Alg.~\ref{alg:proposed}, the following hold with
  $\cost_{n} = \moreau{g}{\mu_{n}}\circ\mathfrak{S}\ (n\in \mathbb{N})$
  and
  $\mu_{n}\coloneqq (2\eta)^{-1}n^{-1/\alpha}$:
  \begin{enumerate}[label=(\alph*)]
    \item
      If
      $\alpha > 1$,
      then
          ${
            \min\limits_{1 \leq n \leq k}
            \mathcal{M}_{\gamma_{n}}^{\cost_{n},\iota_{C}}(\bm{x}_{n})
          \leq \sqrt{\frac{\chi}{(k+1)^{1-\alpha^{-1}} - 1}}}$
          $(k\in\mathbb{N})$
          holds with
      $\chi\coloneqq \frac{\left(\cost(\bm{x}_{1}) - \inf_{\bm{x}\in C} \cost(\bm{x}) + L_{g}^{2}(2\eta)^{-1} \right)(\varpi_{1} + 2\eta\varpi_{2})}{2 \eta c \min\left\{\widetilde{\gamma} (\varpi_{1} + 2\eta\varpi_{2}), \rho(1-2c)\right\}} \in \mathbb{R}_{++}$,
      where
          $\varpi_{1},\varpi_{2} \in \mathbb{R}_{++}$ 
          are
          given
          in Fact~\ref{fact:Lipschitz}.
    \item
          ${\displaystyle\liminf_{n\to\infty} \mathcal{M}_{\gamma_{n}}^{\cost_{n},\iota_{C}}(\bm{x}_{n}) = 0}$.
          \item
          Choose a subsequence
          $(\bm{x}_{m(l)})_{l=1}^{\infty} \subset \mathcal{X}$
          of
          $(\bm{x}_{n})_{n=1}^{\infty}$
          satisfying
          $\lim_{l\to\infty} \mathcal{M}_{\gamma_{m(l)}}^{\cost_{m(l)},\iota_{C}}(\bm{x}_{m(l)}) = 0$.
          Then, every cluster point
          $\bm{x}^{\star} \in C$
          of
          $(\bm{x}_{m(l)})_{l=1}^{\infty}$
          is a stationary point of
          Problem~\ref{problem:origin}.
  \end{enumerate}
\end{theorem}

\begin{remark}[Extension of Alg.~\ref{alg:proposed} to the minimization of $\cost+\phi$ with a prox-regular function $\phi$] \label{remark:extension}
  For simplicity, we are focusing on Problem~\ref{problem:origin} in this paper.
  Nevertheless, Lemma~\ref{lemma:stationarity}, and Thms.~\ref{theorem:asymptotic_approximation} and~\ref{theorem:convergence_extension} can be extended to a case
  where
  $\iota_{C}$
  is replaced by a proper lower semicontinuous and
    {\em prox-regular}\cite[Def. 13.27]{Rockafellar-Wets98} function
  $\phi:\mathcal{X}\to\exR$
  such that
  $\Fsubdiff \phi$
  is {\em outer semicontinuous}~\cite[Def. 5.4]{Rockafellar-Wets98} at every
  $\widebar{\bm{x}} \in \dom{\phi}$
  (Note: such a
  $\phi$
  includes, e.g., the $\ell_{0}$-pseudonorm).
  In this case, the measure
  in~\eqref{eq:measure} is extended to
  $
\thickmuskip=0.2\thickmuskip
\medmuskip=0.2\medmuskip
\thinmuskip=0.2\thinmuskip
\arraycolsep=0.2\arraycolsep
  \mathcal{M}_{\gamma}^{\cost,\phi}(\widebar{\bm{x}}) \coloneqq \inf\left\{\norm{(\widebar{\bm{x}} - \bm{p})/\gamma} \mid \bm{p}\in \prox{\gamma\phi}(\widebar{\bm{x}}-\gamma \bm{v}), \ \bm{v} \in \Fsubdiff \cost(\widebar{\bm{x}}) \right\}$.
  To the problem for finding a stationary point of
  $\cost+\phi$,
  we can apply a modified version of Alg.~\ref{alg:proposed}
  by replacing
  (i)
  $P_{C}$
  in line~\ref{line:projection}
  with
  $\prox{\rho^{m}\widetilde{\gamma}}$;
  and
  (ii)
  the condition in line~\ref{line:Armijo} with the condition
  $(\cost_{n}+\phi)(\bm{x}) \leq (\cost_{n}+\phi)(\bm{x}_{n}) - c\rho^{m}\widetilde{\gamma} \norm{( \bm{x}_{n}-\bm{x} )/(\rho^{m}\widetilde{\gamma})}^{2}$.
\end{remark}

\section{Numerical Experiments on Synthetic Data}

In a scenario of RLRMR in Problem~\ref{problem:RLRMR},
we demonstrate the performance of the proposed model~\eqref{eq:proposed_model} with
$\sigma \coloneqq 1$
for
$\mathcal{L}_{r,\sigma}$
in~\eqref{eq:low_rank}
and with two loss functions:
(i) the $\ell_{1}$-norm
$\norm{\cdot}_{1}$,
i.e.,
$\ell \coloneqq \abs{\cdot}$,
and
(ii) SCAD~\cite{Fan-Li01}, i.e.,
$\ell\coloneqq r^{\rm SCAD}_{\theta}$
with a parameter
$\theta > 2$:
\begin{equation}
  r^{\rm SCAD}_{\theta}:\mathbb{R}\to\mathbb{R}:t\mapsto 
  \begin{cases}
    \abs{t},                                     & \mathrm{if}\ \abs{t}\leq 1;                    \\
    \frac{-t^{2}+2\theta\abs{t}-1}{2(\theta-1)}, & \mathrm{if}\ 1< \abs{t} \leq \theta;           \\
    \frac{\theta+1}{2},                          & \mathrm{if}\ \abs{t} > \theta, \label{eq:SCAD}
  \end{cases}
\end{equation}
which is $(\theta-1)^{-1}$-weakly convex (e.g.,~\cite{Bohm-Wright21}).
For comparisons, we employed two state-of-the-art models: the model~\eqref{eq:RLRMR_decomp} in~\cite{Li-Zhu-Man-Vidal20} and the following model in~\cite{Xu-Li-Zheng24}
with weights
$\lambda,\beta \in \mathbb{R}_{++}$:
\begin{equation}
  \mathop{\mathrm{minimize}}_{\bm{X}\in\mathbb{R}^{n_{1}\times n_{2}}} \ \norm{\bm{y}-\mathcal{A}(\bm{X})}_{1} +\lambda (\norm{\bm{X}}_{\rm nuc} - \beta \norm{\bm{X}}_{F}), \label{eq:NRLAD}
\end{equation}
where
$\norm{\cdot}_{\rm nuc}$
denotes the nuclear norm, and
the second term in~\eqref{eq:NRLAD} is
a regularizer to promote the low-rankness of
$\bm{X}$.
We applied (a) Alg.~\ref{alg:proposed} with
$(c,\rho,\widetilde{\gamma},\alpha) \coloneqq (2^{-13},2^{-1},1,3)$
to~\eqref{eq:proposed_model}
(see also Remark~\ref{remark:formulation});
(b) a subgradient method~\cite[Alg. 2.1]{Li-Zhu-Man-Vidal20} (see also~\cite{Davis-Drusvyatskiy-MacPhee18}) to~\eqref{eq:RLRMR_decomp};
and
(c)
a difference-of-convex algorithm (DCA)~\cite[Alg. 1]{Xu-Li-Zheng24} (see also~\cite{Le-Pham24}) to~\eqref{eq:NRLAD}, where an ADMM, e.g.,~\cite{He-Yuan12}, was employed for a convex subproblem in DCA~\cite[Alg. 1]{Xu-Li-Zheng24}.
All algorithms were terminated when runtime exceeded $60$ (s) or the relative difference
$\abs{J(\bm{x}_{n+1})-J(\bm{x}_{n})}/\abs{J(\bm{x}_{n})}$
was less than
$10^{-9}$,
where
$J$
is the cost function in~\eqref{eq:RLRMR_decomp},~\eqref{eq:proposed_model}, or~\eqref{eq:NRLAD},
and
$\bm{x}_{n}$
is an
$n$th estimate of a solution generated by each algorithm.
All experiments were performed by MATLAB on MacBookPro (Apple M3, 16 GB).

For every trial, according to~\cite{Li-Zhu-Man-Vidal20,Xu-Li-Zheng24}, we generated
$\bm{X}^{\star} \coloneqq \bm{U}^{\star}\bm{V}^{\star\TT}\in \mathcal{L}_{r} \subset \mathbb{R}^{n_{1}\times n_{2}}$
with random
$\bm{U}^{\star} \in \mathbb{R}^{n_{1}\times r}$
and
$\bm{V}^{\star} \in \mathbb{R}^{n_{2}\times r}$,
$\mathcal{A}_{i}:\mathbb{R}^{n_{1}\times n_{2}} \to \mathbb{R}\ (i=1,2,\ldots,m)$,
and
$\bm{y} \in \mathbb{R}^{m}$
as follows.
Each entry of
$\bm{U}^{\star}$,
$\bm{V}^{\star}$
and all
$\mathcal{A}_{i}$,
as matrix expressions of
$\mathcal{A}_{i}$,
was sampled from the normal distribution
$\mathcal{N}(0,1)$.
We generated
$\bm{y} \in \mathbb{R}^{m}$
by~\eqref{eq:observation},
where each entry of noise
$\bm{\epsilon} \in \mathbb{R}^{m}$
was sampled from
$\mathcal{N}(0,10^{-6})$
according to~\cite{Xu-Li-Zheng24},
the index set
$\mathcal{I}_{\rm out}$
of outliers was randomly chosen, and each outlier
$\xi_{i} \in \mathbb{R}\ (i\in\mathcal{I}_{\rm out})$
was generated with
$\Omega \coloneqq \max_{i=1,2,\ldots,m}{\abs{\mathcal{A}_{i}(\bm{X}^{\star})}}$
from
(i) uniform distribution
$\mathcal{U}_{[-\Omega, \Omega]}$;
and
(ii) Cauchy distribution, i.e.,
$\xi_{i}\coloneqq \Omega \tan(\pi u_{i}/2)$
with
$u_{i} \in (-1,1)$
sampled from
$\mathcal{U}_{(-1,1)}$.

Table~\ref{table:result} demonstrates the averaged root-mean-square error (RMSE) and runtime of each model\footnote{
For~\eqref{eq:proposed_model} with
$\ell\coloneqq r_{\theta}^{\rm SCAD}$ in~\eqref{eq:SCAD}, we employed
$\theta \in \{2.5,2.7,2.9,3.1\}$.
For~\eqref{eq:RLRMR_decomp},
$\lambda \in \{10^{-2},10^{-1},1, 2, 5\}$
was chosen.
Moreover, stepsize of the subgradient method used in~\cite{Li-Zhu-Man-Vidal20} at $n$th iteration was chosen from
$0.95^{n}$
or
$2\times 0.95^{n}$,
where the first choice was suggested by~\cite{Li-Zhu-Man-Vidal20} while the second choice achieved better results in some cases in our scenario.
For~\eqref{eq:NRLAD},
$\lambda \coloneqq t \sqrt{mn_{2}\log(n_{1}+n_{2})}$
with
$t \in \{0.1, 0,5\}$
and
$\beta \in \{0.1, 0.5, 0.9\}$
were chosen according to~\cite{Xu-Li-Zheng24}.
} over $100$ trials with
$(n_{1},n_{2},r) = (40,50,5)$,
$m = p_{m}n_{1}n_{2}$
and
$\abs{\mathcal{I}_{\rm out}} = p_{\abs{\mathcal{I}_{\rm out}}}m$
for
$p_{m} \in \{0.3,0.5,0.7\}$
and
$p_{\abs{\mathcal{I}_{\rm out}}} \in \{0.7,0.5\}$,
where we used RMSE
$\lVert\widehat{\bm{X}}-\bm{X}^{\star}\rVert_{F}/\sqrt{n_{1}n_{2}}$
to measure a recovery error with the final estimate
$\widehat{\bm{X}} \in \mathbb{R}^{n_{1}\times n_{2}}$
of algorithm.
We note that Problem~\ref{problem:RLRMR} becomes challenging as
$p_{m}$
is smaller and
$p_{\abs{\mathcal{I}_{\rm out}}}$
is larger.

From Table~\ref{table:result}, we observe that convergence speeds of the proposed models~\eqref{eq:proposed_model} and the model~\eqref{eq:NRLAD} are much slower than that of the model~\eqref{eq:RLRMR_decomp}.
This is because singular value decompositions (SVD) of
$n_{1}$-by-$n_{2}$
matrices
are required at every iteration in solvers for~\eqref{eq:proposed_model} (see also Remark~\ref{remark:formulation}) and for~\eqref{eq:NRLAD} while
any SVD is not required in a subgradient method for~\eqref{eq:RLRMR_decomp} because of the decomposition
$\bm{X} = \bm{U}\bm{V}^{\TT}$
in~\eqref{eq:RLRMR_decomp}
(Note: acceleration of the proposed Alg.~\ref{alg:proposed} is beyond the scope of this paper, and will be addressed in future work).
Moreover, from Table~\ref{table:result},
the proposed model~\eqref{eq:proposed_model} with
$\ell\coloneqq r^{\rm SCAD}_{\theta}$
achieves smaller RMSE than the model~\eqref{eq:proposed_model} with
$\ell\coloneqq \abs{\cdot}$
in a shorter runtime;
and
RMSE of the existing models~\eqref{eq:RLRMR_decomp} and~\eqref{eq:NRLAD} are much higher than that of the model~\eqref{eq:proposed_model} with
$\ell\coloneqq r^{\rm SCAD}_{\theta}$
in particular for
$p_{m} \in \{0.3,0.5\}$.
These observations imply that
a weakly convex
$r^{\rm SCAD}_{\theta}$
is more robust against outliers,
at least in the proposed model~\eqref{eq:proposed_model},
than a convex
$\abs{\cdot}$
as we expected; and
the proposed model~\eqref{eq:proposed_model} with
$\ell\coloneqq r^{\rm SCAD}_{\theta}$
outperforms the existing models~\eqref{eq:RLRMR_decomp} and~\eqref{eq:NRLAD} in terms of recovery accuracy\footnote{
  We note that theoretical recovery guarantees of optimization models involving nonconvex loss functions, e.g., SCAD, have been reported for special cases of RLRMR~\eqref{eq:RLRMR},
  e.g.,
  robust matrix completion~\cite{Wang-Wei24}.
  However, such a theoretical recovery guarantee of the proposed model~\eqref{eq:proposed_model} is beyond the scope of this paper, and will be addressed in future work.
}.

  {
    \tabcolsep = 0.2em
    \begin{table}[t]
      \caption{\footnotesize Averaged RMSE and (averaged runtime [s]) for Problem~\ref{problem:RLRMR}
      }
      \label{table:result}
      \hspace{-1.4em}
      {
        \scriptsize
        \begin{tabular}{lr|ll|ll}\toprule
          \multicolumn{2}{c|}{Setting}
           & \multicolumn{2}{c|}{Proposed model~\eqref{eq:proposed_model}}

           & \multicolumn{2}{c}{Existing models}                           \\
          
          $(p_{m},p_{\abs{\mathcal{I}_{\rm out}}})$
           & outlier
           & $\ell\coloneqq \abs{\cdot}$

           & $\ell\coloneqq r_{\theta}^{\rm SCAD}$

           & Model~\eqref{eq:RLRMR_decomp}~\cite{Li-Zhu-Man-Vidal20}

           & Model~\eqref{eq:NRLAD}~\cite{Xu-Li-Zheng24}                   \\ \midrule
          
          \multirow{2}{*}{$(0.3, 0.7)$}
           & uniform
           & 6.81E-01 {\notsotiny (59.9)}
           & {\bf 3.46E-05} {\notsotiny (21.8)}

           & 9.79E-01 {\notsotiny (0.7)}
           & 1.89E-01 {\notsotiny (53.3)}                                  \\
           & Cauchy 
           & 6.56E-01 {\notsotiny (59.8)}
           & {\bf 3.47E-05} {\notsotiny (23.6)}

           & 9.93E-01 {\notsotiny (0.6)}
           & 4.59E-01 {\notsotiny (58.7) }                                 \\ \midrule

          \multirow{2}{*}{$(0.3, 0.5)$}
           & uniform
           & 7.02E-01 {\notsotiny (60.0)}
           & {\bf 3.83E-05} {\notsotiny (20.2)}

           & 1.04E+00 {\notsotiny (0.7)}
           & 1.83E-01 {\notsotiny (55.9) }                                 \\
           & Cauchy 
           & 6.87E-01 {\notsotiny (60.0)}
           & {\bf 3.81E-05} {\notsotiny (20.3)}

           & 1.10E+00 {\notsotiny (0.7)}
           & 1.69E-01 {\notsotiny (55.8)}                                  \\ \midrule
          
          \multirow{2}{*}{$(0.5, 0.7)$}
           & uniform
           & 5.09E-05 {\notsotiny (58.4)}
           & {\bf 1.92E-05} {\notsotiny (5.5)}

           & 2.92E-02 {\notsotiny (1.1)}
           & 4.78E-03 {\notsotiny (59.3)}                                  \\
           & Cauchy 
           & 5.47E-05 {\notsotiny (50.9)}
           & {\bf 1.91E-05} {\notsotiny (6.9)}

           & 7.39E-03 {\notsotiny (1.2)}
           & 4.53E-03 {\notsotiny (58.5)}                                  \\ \midrule
          
          \multirow{2}{*}{$(0.5, 0.5)$}
           & uniform
           & 5.02E-05 {\notsotiny (57.7)}
           & {\bf 1.92E-05} {\notsotiny (6.2)}

           & 8.11E-03 {\notsotiny (1.0)}
           & 5.52E-03 {\notsotiny (59.3)}                                  \\
           & Cauchy 
           & 5.48E-05 {\notsotiny (50.7)}
           & {\bf 1.98E-05} {\notsotiny (5.7)}

           & 3.61E-03 {\notsotiny (1.1)}
           & 5.05E-03 {\notsotiny (58.0)}                                  \\ \midrule
          
          \multirow{2}{*}{$(0.7, 0.7)$}
           & uniform
           & 3.17E-05 {\notsotiny (58.6)}
           & {\bf 1.50E-05} {\notsotiny (4.4)}

           & 1.72E-05 {\notsotiny (1.3)}
           & 1.35E-03 {\notsotiny (7.4)}                                   \\
           & Cauchy 
           & 3.34E-05 {\notsotiny (52.8)}
           & {\bf 1.48E-05} {\notsotiny (4.0)}

           & 1.74E-05 {\notsotiny (1.6)}
           & 1.88E-03 {\notsotiny (9.2)}                                   \\ \midrule
          
          \multirow{2}{*}{$(0.7, 0.5)$}
           & uniform
           & 3.15E-05 {\notsotiny (58.2)}
           & {\bf 1.48E-05} {\notsotiny (4.0)}

           & 1.75E-05 {\notsotiny (1.7)}
           & 1.31E-03 {\notsotiny (8.4)}                                   \\
           & Cauchy 
           & 3.18E-05 {\notsotiny (55.9)}
           & {\bf 1.51E-05} {\notsotiny (5.1)}

           & 1.74E-05 {\notsotiny (1.4)}
           & 1.97E-03 {\notsotiny (9.1)}                                   \\ \bottomrule
          
        \end{tabular}
      }
    \end{table}
  }

\section{Concluding Remarks}
In this paper, we proposed a prox-regular-type low-rank constrained optimization model, incorporating a weakly convex loss function, for robust low-rank matrix recovery.
For the proposed model, we presented an optimization algorithm, for minimization of a weakly convex function over a prox-regular set, with guaranteed convergence in terms of stationary point.
Our numerical experiments demonstrate that the proposed model with a weakly convex loss function, called SCAD, dramatically improves estimation performance of robust low-rank matrix recovery under severe outliers.
Finally, we defer further numerical evaluations and complete proofs to an extended manuscript in preparation.

\bibliographystyle{IEEEtran}
\bibliography{main}

\end{document}